\documentclass{amsart}
\usepackage{amsmath}
\usepackage{amsfonts}
\usepackage{mathrsfs}
\usepackage{amssymb}
\usepackage{amsthm}
\newtheorem{theorem}{Theorem}[section]
\newtheorem{lemma}[theorem]{Lemma}
\newtheorem{conjecture}[theorem]{Conjecture}

\theoremstyle{definition}

\theoremstyle{remark}
\newtheorem{remark}[theorem]{Remark}

\numberwithin{equation}{section}

\newcommand{\BBR}{{\mathbb R}}

\newcommand{\BBN}{{\mathbb N}}
\newcommand{\BBC}{{\mathbb C}}
\newcommand{\BBZ}{{\mathbb Z}}
\newcommand{\BBT}{{\mathbb T}}
\newcommand{\T}{{\operatorname T}}
\newcommand{\M}{{\operatorname M}}
\newcommand{\Ess}{{\operatorname S}}
\newcommand{\E}{{\operatorname E}}

\newcommand{\proj}{{\operatorname{proj}}}
\newcommand{\tr}{{\operatorname{tr}}}
\newcommand{\AM}{{\operatorname{H}_{\lambda, \alpha, \theta}}}
\newcommand{\ch}{{\operatorname{ch}}}

\begin{document}

\title[Linear independence of lattice Gabor systems]{On the finite linear independence of lattice Gabor systems}

\author[C.\ Demeter]{Ciprian Demeter}
\address{Department of Mathematics, Indiana University, Bloomington, IN 47405}
\email{demeterc@indiana.edu}
\thanks{The first author is supported by a Sloan Research Fellowship and by
NSF Grant DMS-0901208.}

\author[S.\ Z.\ Gautam]{S.\ Zubin Gautam}
\address{Department of Mathematics, Indiana University, Bloomington, IN 47405}
\email{sgautam@indiana.edu}



\date{}

\dedicatory{}

\begin{abstract}
In the restricted setting of product phase space lattices, we give
an alternate proof of P.\ Linnell's theorem on the finite linear
independence of lattice Gabor systems in $L^2(\BBR^d)$.  Our proof
is based on a simple argument from the spectral theory of random
Schr\"odinger operators; in the one-dimensional setting, we recover
the full strength of Linnell's result for general lattices.
\end{abstract}

\maketitle

\section{Introduction}
A \emph{Gabor system} is simply a collection of modulations and
translations of a fixed function in $L^2(\BBR^d)$.  More precisely,
given any set $A \subseteq \BBR^d \times \BBR^d$ and any $f \in
L^2(\BBR^d)$, the associated Gabor system is \[\mathcal G(f, A) :=
\{\M_y \T_x f \: | \: (x,y) \in A \subset \BBR^d \times \BBR^d\},\]
where $\M_y$ and $\T_x$ denote, respectively, the unitary operators
of modulation and translation on $L^2(\BBR^d)$ given by \[\M_y f(t)
:= e^{2\pi i y \cdot t}f(t), \:\: \T_x f(t) := f(t-x).\]  Here one
should view $\BBR^d \times \BBR^d \cong \BBR^d \times
\widehat{\BBR^d}$ as the phase space of $\BBR^d$; accordingly, due
to the Fourier transform's intertwining of modulation and
translation, one may view $\mathcal G(f,A)$ as the collection of
``phase space translates of $f$ by $A$'' or, inspired by the case
$d=1$, the collection of ``time-frequency translates of $f$ by
$A$.'' Most of the interest in Gabor systems stems from their
``basis-like'' utility in providing \emph{expansions} of $L^2$
functions; thus, quite naturally, a significant portion of research
in the field has focused on investigating completeness and
independence properties of these systems. Namely, given input
data $(f, A)$ in some particular class, one might like to know
whether one can deduce that $\mathcal G(f,A)$ is an orthonormal
basis, a Schauder basis, a frame, or one of a number of other
basis-like objects for $L^2(\BBR^d)$ (or more generally for the span
of $\mathcal G(f,A)$).

Perhaps the most basic independence property one could ask of
$\mathcal G(f,A)$ is that it be \emph{finitely linearly
independent}.  The \emph{Heil--Ramanathan--Topiwala Conjecture}
(henceforth referred to as the ``HRT Conjecture'') asserts that
\emph{any} nontrivial Gabor system should have this property:

\begin{conjecture}[Heil--Ramanathan--Topiwala]\label{hrtconjecture}  Let $A \subset
\BBR^d \times \BBR^d$ be a finite set and $0\neq f \in
L^2(\BBR^d)$.  Then the associated Gabor system $\mathcal G(f,A)$ is
linearly independent as a subset of $L^2(\BBR^d)$.
\end{conjecture}

This conjecture was originally posed in \cite{hrtpaper}, in which
the claim was verified under various restrictions on either the set
$A$ or the function $f$ (see also \cite{heil} for a nice expository account); two of these results in particular point
the way to the setting of this paper, \textit{viz.\ }that of
\emph{lattice} Gabor systems.  Namely, in the one-dimensional $d=1$
setting, Heil, Ramanathan, and Topiwala proved Conjecture
\ref{hrtconjecture} under the restriction that $A$ be an arbitrary
finite subset of a covolume-$1$ lattice in $\BBR \times \BBR$, as
well as under the alternate restriction that $A$ have cardinality at
most $3$, which in turn implies that $A$ or some translate thereof
is contained in a lattice in $\BBR \times \BBR = \BBR^2$.  (We
recall that a \emph{lattice} in $\BBR^n$ is a discrete subgroup
$\Gamma \leq \BBR^n$ of finite covolume; we define the covolume of
$\Gamma$ to be the Lebesgue measure of a fundamental domain for the
quotient space $\BBR^n / \Gamma$.)

The HRT Conjecture is of course strikingly simple in formulation;
however, after the preliminary results proved in \cite{hrtpaper},
surprisingly little headway has been made toward its resolution.
Indeed, the only landmark result of a reasonably general nature
regarding this conjecture is the following 1999 theorem of Linnell,
which carries the aforementioned results of \cite{hrtpaper} to a
natural conclusion:
\begin{theorem}[Linnell (\cite{linnell99})]\label{linnelltheorem}
Suppose $0 \neq f \in L^2(\BBR^d)$, and suppose that some translate
of a finite set $A \subset \BBR^d \times \BBR^d$ is contained in a
lattice of $\BBR^d \times \BBR^d$.  Then the Gabor system $\mathcal
G(f,A)$ is linearly independent in $L^2(\BBR^d)$.
\end{theorem}

Linnell's proof is based on a ``twisted'' version of the group von
Neumann algebra techniques he developed in \cite{linnell91} to prove
the (analytic) zero divisor conjecture for elementary amenable
groups. While the particular von Neumann algebras exploited by
Linnell arise rather naturally in Gabor analysis,
there has been interest in obtaining a more ostensibly
``elementary'' proof of Theorem \ref{linnelltheorem}, toward a
better understanding of the HRT Conjecture.  An alternate proof of
the \emph{one-dimensional} ($d=1$) case was recently given by Bownik
and Speegle (\cite{bownikspeegle}) using the theory of
shift-invariant spaces; however, as noted by the authors, the
methods do not extend readily to higher dimensions.

The purpose of this paper is to provide yet another proof of Theorem
\ref{linnelltheorem} that is valid for arbitrary dimensions $d$;
unfortunately, however, our proof seems only to be able to treat the
case of \emph{product lattices} of the form $\Gamma \times \Lambda
\leq \BBR^d \times \BBR^d$, or more generally lattices that can be
mapped to a product lattice by a symplectic transformation of the
phase space.  Nonetheless, we note that this case is \emph{generic}
in the $d=1$ setting:  \emph{any} lattice in $\BBR \times \BBR$ can
be symplectically mapped to a product lattice; we will discuss these
issues in greater detail in Section \ref{metaplectic} below.  Our
main result is the following:

\begin{theorem}\label{maintheorem}
Let $\Gamma_0 = \Gamma \times \Lambda$ be a \emph{product} lattice
in $\BBR^d \times \BBR^d$, where $\Gamma$ and $\Lambda$ are
arbitrary lattices in $\BBR^d$.  Suppose that some translate of a
finite set $A \subset \BBR^d \times \BBR^d$ is contained in
$\Gamma_0$, and suppose $0 \neq f \in L^2(\BBR^d)$. Then the Gabor
system $\mathcal G(f,A)$ is linearly independent.  The same result
holds if $\Gamma_0$ is merely assumed to be a symplectic image of a
product lattice; in particular, for $d=1$ the result holds for
$\Gamma_0$ an arbitrary lattice in $\BBR \times \BBR$.
\end{theorem}

Because of the product lattice restriction in our theorem, Linnell's
proof remains the state of the art regarding the HRT conjecture for
lattices; however, our proof does have the advantage of actually
yielding a slightly stronger result in this restricted setting.  To
wit, once one has reduced $\Gamma_0$ to a product lattice by a
symplectic transformation, the particular structure of the
modulation operators in the definition of a Gabor system becomes
inconsequential; namely, one can replace the modulations with
multiplication operators by \emph{any}
almost-everywhere-nonvanishing $L^\infty$ functions with a suitable
periodicity.\footnote{Given Linnell's proof this strengthening is
intuitively not so surprising, since the more general operators that
arise are obviously contained in the von Neumann algebra generated
by $\{\M_y \T_x \; | \; (x,y) \in \Gamma_0\}$. However, it is not
immediately clear to us whether or how easily the full strength of
Theorem \ref{generaltheorem} below can be obtained by the methods of
\cite{linnell99}; in short, some of the ``Ore localization'' issues
arising in Linnell's proof seem to become less trivial after passing from
modulations to more general multiplication operators.} See Theorem
\ref{generaltheorem} below for a precise statement of this
generalization.  Furthermore, one can of course consider Gabor systems $\mathcal G(f,A)$ whose generating functions $f$ are not elements of $L^2$; at the expense of sacrificing symplectic symmetry, our method of proof immediately applies to a larger class of such Gabor systems:

\begin{theorem}\label{lptheorem}
Suppose $0<p\leq 2$.  If some translate of  $A \subset \BBR^d \times \BBR^d$ is contained in a \emph{product} lattice, then the Gabor system $\mathcal G(f,A) \subset L^p(\BBR^d)$ is linearly independent for any $0 \neq f \in L^p(\BBR^d)$.
\end{theorem}

\smallskip
The methods of this paper were inspired by a connection between the
HRT Conjecture and discrete Schr\"odinger operators that was
observed by F.\ Nazarov and A.\ Volberg and pointed out to us by C.\
Thiele. Specifically, for parameters $\lambda > 0$ and $\alpha,
\theta \in \BBT = \BBR / \BBZ \cong \widehat \BBZ$, consider the
\emph{almost Mathieu operator} $\operatorname{H}_{\lambda, \alpha,
\theta}$ on $\ell^2 \BBZ$ defined by
\begin{align*}
\AM u(n) &= u(n+1) + u(n-1) + 2\lambda \cos 2\pi (\theta + n \alpha)
\cdot u(n)\\
&= u(n+1) + u(n-1) + [\lambda e^{2\pi i \theta}] e^{2\pi i \alpha n}
u(n) + [\lambda e^{-2\pi i \theta}] e^{2\pi i (-\alpha) n} u(n).
\end{align*}
Observe that if $0\neq u \in \ell^2 \BBZ$ satisfies an eigenvalue
equation
\begin{equation}\label{eigenvalue}
\AM u = Eu,
\end{equation}
then one can consider the set $A = \{(-1,0), (1,0),
(0,\alpha), (0, - \alpha), (0,0)\} \subset \BBZ \times \widehat
\BBZ$ to obtain a linearly dependent analogue $\mathcal G(u,A)$ of a
Gabor system, with the group $\BBR^d$ replaced by $\BBZ$.  (We
define $\mathcal G(u,A)$ in the obvious way via translations and
modulations of $u$, with multiplication by characters of $\BBZ$
giving the modulations.)  That is, the ``Gabor system'' $\mathcal
G(u,A)$ gives a \emph{counterexample} to the analogue of the HRT
Conjecture over the group $\BBZ$.

The spectral theory of the almost Mathieu operator has long been a
focal point in the study of discrete Schr\"odinger operators. For
our present purposes, it suffices merely to note that for certain
values of $\lambda, \alpha$, and $\theta$ the operator $\AM$ has
some pure-point spectrum; indeed, for $\lambda > 1$ the spectrum was
shown by Jitomirskaya to be \emph{entirely} pure-point for almost
every $\alpha, \theta \in \BBT$ (see \cite{jitomirskaya}).  That is,
one can actually find nonzero eigenfunctions $u$ and obtain
counterexamples to the HRT Conjecture over $\BBZ$.  On the other
hand, one could recast the set $A$ above as a subset of the
time-frequency lattice $\BBZ \times \alpha \BBZ \leq \BBR \times
\widehat \BBR$ and consider a Gabor system $\mathcal G(f,A)$ with
$0\neq f \in L^2(\BBR)$; this Gabor system must be linearly
independent by Theorem \ref{linnelltheorem}.\footnote{We must admit
that this recasting is a bit artificial; $A$ is a subset of a
\emph{lattice} in the phase space over $\BBZ$ precisely when
$\alpha$ is \emph{rational}, so for $\alpha$ irrational $A$ takes on
a decidedly different nature when viewed as a subset of $\BBR \times
\widehat \BBR$.  For $\alpha$ rational the methods below easily
yield the absence of pure-point spectrum for $\AM$ (in fact the
spectrum is well known to be entirely absolutely continuous in that
case), so $\AM$ itself does not yield a counterexample to the HRT
Conjecture over $\BBZ$.  We also note that failure of the full HRT
Conjecture over $\BBZ$ can be deduced by much simpler methods; indeed, taking $A' = \{(0, \alpha), (0, \beta)\}$ for any $\alpha \neq \beta \in \widehat \BBZ$, the Gabor system $\mathcal G(\delta_n, A')$ is trivially linearly dependent for any $0 \neq n \in \BBZ$.  Cf.\ Theorem 1 of \cite{kutyniok} for a more general discussion.}

This contrast motivates our approach below,
which essentially begins by decomposing the function $f$ along
$\BBZ$-orbits $\{f(x+n)\}_{n\in \BBZ}$ and examining hypothetical
linear dependence relations inside the resulting discrete Gabor
systems; we show that (essentially) no such relations can occur.  Indeed, in
the context of the almost Mathieu operators, suppose that for fixed $\lambda$ and $\alpha$ one could find eigenfunctions $u_\theta \in \ell^2\BBZ$ of $\AM$ satisfying (\ref{eigenvalue}) for some eigenvalue $E$ \emph{independent} of almost every $\theta \in \BBT$.  Then one could ``piece together'' the $u_\theta$ appropriately to obtain $f\in L^2(\BBR)$ with $\mathcal G(f,A)$ violating Linnell's Theorem \ref{linnelltheorem}, so the eigenvalues $E$ must actually depend on $\theta$.  In fact, the dependence must be quite ``wild'' in the sense that any fixed $E \in \BBR$ is $\theta$-almost surely \emph{not} an eigenvalue of $\AM$; this fact was already known prior to Linnell's result, and the most of the main arguments of Section \ref{proofsection} below closely follow its standard proof (cf.\ \textit{e.g.\ }Proposition V.2.8, \cite{carmonalacroix}).  Thus, the most basic
perspective from which one should approach our proof is the
following: an operator from the algebra yielding a product lattice
Gabor system naturally gives rise to a measurable family of
operators that is rather similar to a random discrete Schr\"odinger
operator.  Moreover, the study of recurrence relations arising from a discrete translation group structure, which is central to our argument, has already yielded new results in the setting of non-lattice Gabor systems; see \cite{demetermn} and \cite{demeterzaharescu}.

\smallskip
At this point, we should briefly discuss the common spirit of our
proof for product lattices, Linnell's proof for general lattices,
and the Bownik--Speegle proof of the $d=1$ case.  Indeed, there is
considerable aesthetic similarity between our argument and that of
Bownik and Speegle, and readers familiar with the basics of von
Neumann algebras will notice such a similarity among all three
proofs.  The core phenomenon behind all the arguments is the
following:  On the one hand, one uses basic generalities to obtain
\textit{a priori} estimates on the dimensions of certain spaces
(spans of Gabor systems in the Bownik--Speegle setting, and kernels
of operators in ours and Linnell's); on the other hand, one shows
that a linear dependence violating the HRT Conjecture would force
the appropriate dimensions to be either smaller or larger than
reality permits.  The current paper and \cite{bownikspeegle} use the
classical notion of dimension of subspaces of $L^2$, while Linnell
considers the Murray--von Neumann dimension of subspaces (\textit{i.e.}, the
Murray--von Neumann trace of their associated orthogonal projections) relative
to the aforementioned von Neumann algebras.

In the exact same spirit, we should also point out an elegant
unpublished proof due to Thiele (\cite{thielecomm}), which yields
Theorem \ref{linnelltheorem} under sufficient time-frequency decay
conditions on $f$ (demanding that $f$ be Schwartz is more than
sufficient).  Thiele's argument takes any lattice $\Gamma_0 \leq
\BBR^d \times \BBR^d$, fixes a function $f \in L^2(\BBR^d)$, and
considers the Gabor systems associated to $f$ and large balls of
radius $R$ in the lattice $\Gamma_0$.  If one supposes a linear
dependence in $\mathcal G(f,A)$ for some $A \subset \Gamma_0$, then
the property of free abelian groups appearing in Lemma
\ref{propagationlemma} below, which we use in a different context, shows that
the spans of these Gabor systems have dimension $O(R^{2d-1})$.  On
the other hand, via an almost-orthogonality argument, the time-frequency
decay assumptions on $f$ show that the dimensions must actually grow
more quickly than $cR^{2d-1 + \varepsilon}$.

\smallskip
We conclude this introduction by pointing out one final advantage of
our proof, namely that it showcases a hierarchy of complexity in the
HRT problem for lattices.  Firstly, there is a jump in complexity
from the one-dimensional case to that of higher dimensions, not only
in the facility of metaplectic reduction (cf.\ Section
\ref{metaplectic}) but also in obtaining \emph{a priori}
bounded-dimensionality conditions.\footnote{This increase in
complexity is already manifest for different reasons in Linnell's
argument, which can be made considerably more ``concrete'' for $d=1$;
much of the subtlety of his proof lies in the induction step
required to extend past this case. However, given that his result is
often cited only for $d=1$, this feature may have been somewhat
underappreciated in the literature.} Secondly, in the
higher-dimensional setting, there is likewise a jump in complexity
between the product lattice setting and that of general lattices;
this is most obviously illustrated by the failure of our arguments
outside the product setting. Perhaps more interesting, however, is
the fact that once one has ``decoupled'' a lattice into a product
lattice in phase space (if possible), anything resembling Fourier
analysis or the structure of the Heisenberg group (cf.\ Section \ref{metaplectic}) disappears from
the picture. This is apparent in the aforementioned fact that one
need not consider actual \emph{modulations} in the product setting; see also Remark \ref{groupremark} of Section \ref{remarks} for a further discussion in this vein.

In the following, for two quantities $A$ and $B$, we will use the notation ``$A \lesssim B$'' to denote the inequality $A \leq cB$ for some constant $c$.  Whenever necessary, any dependence of the implied constants $c$ on relevant parameters will be denoted by subscripts on the symbol ``$\lesssim$.''

\subsection*{Acknowledgements}  We would like to thank Nets Katz and
Christoph Thiele for many useful discussions regarding the HRT
Conjecture.  We are especially grateful to Christoph Thiele both for
pointing out the connection with the almost Mathieu operator and for
showing us the aforementioned proof for functions with
time-frequency decay; the current work would have been impossible
without either of these observations.


\section{Metaplectic reductions}\label{metaplectic}
A crucial observation appearing in \cite{hrtpaper} is that the
linear independence of Gabor systems is invariant under certain
affine transformations of the phase space; more precisely, suppose
that for some $A \subset \BBR^{2d} \cong \BBR^d \times
\widehat{\BBR^d}$ we know $\mathcal G(f,A)$ is linearly independent
for all nonzero $f\in L^2(\BBR^d)$. Then we automatically know that
$\mathcal G(f,\sigma A)$ is independent for all $f$ whenever $\sigma
\in \operatorname{Sp}_{2d}(\BBR) \ltimes \BBR^{2d}$ is an
\emph{affine-symplectic} transformation, \textit{i.e.\ }$\sigma$ is
a is a composition of a translation and a linear transformation
preserving the symplectic form on $\BBR^{2d}$. This symmetry is
essentially due to the fact that replacing $A$ by $\sigma A$ amounts
to pre- and post-composing the operators $\M_y \T_x$, $(x,y) \in A$, by some
unitary ``metaplectic transformations'' in $\mathcal U\big(L^2(\BBR^d)\big)$; this in turn is due to the fact that
the linear action $\operatorname{Sp}_{2d}(\BBR) \curvearrowright
\BBR^{2d}$ induces automorphisms of the Heisenberg group with
underlying set $\BBR^{2d} \times \BBR$, and Gabor systems in
$L^2(\BBR^d)$ arise from the (unitary) Schr\"odinger representation
of this Heisenberg group. For a more complete discussion, see
\textit{e.g.\ }Section XII.7.B of \cite{bigstein}.  The translation
symmetry in particular accounts for the ``some translate of $A$''
phrasing in Theorems \ref{linnelltheorem} and \ref{maintheorem};
henceforth we will restrict our attention to actual subsets of
lattices.

Accordingly, we declare two lattices $\Gamma_1, \Gamma_2 \leq
\BBR^{2d}$ to be \emph{symplectically equivalent} if there is some
$\sigma \in \operatorname{Sp}_{2d}(\BBR)$ for which $\sigma \Gamma_1
= \Gamma_2$.  In the case $d=1$, the symplectic group
$\operatorname{Sp}_2(\BBR)$ luckily coincides with the entire
special linear group $\operatorname{SL}_2(\BBR)$, which is easily
seen to act \emph{transitively} on the space of lattices of a given
covolume in $\BBR^2$.  Thus, in particular, \emph{any lattice
$\Gamma_1 \leq \BBR^2$ is symplectically equivalent to a product
lattice $\Gamma_2 \leq \BBR \times \BBR$}.

For $d \geq 2$, however, $\operatorname{Sp}_{2d}$ is a \emph{proper}
subgroup of $\operatorname{SL}_{2d}$, and in fact one can construct
lattices in $\BBR^{2d}$ that are \emph{not} symplectically
equivalent to any product lattice (cf.\ Remark \ref{productremark} below).  Thus,
we resign ourselves to the restricted setting of product lattices,
and Theorem \ref{maintheorem} only treats generic lattices for
$d=1$.

\section{Proof of the main theorem}\label{proofsection}
At last, we come to the proof of Theorem \ref{maintheorem}; by the
discussion of the previous section, we need only treat the case in
which $\Gamma_0$ is a genuine product lattice.  Our most basic
perspective is identical to that of \cite{linnell99}; namely, in
lieu of studying the Gabor systems $\mathcal G(f,A)$ themselves, we
examine the algebra of operators generated by $\{\M_y \T_x \, | \,
(x,y) \in \Gamma_0\}$ in the space $\mathcal B\big(L^2(\BBR^d)\big)$
of bounded operators on the Hilbert space $L^2(\BBR^d)$.  Indeed,
suppose $f \in L^2(\BBR^d)$ is in the \emph{kernel} of some operator
$\Ess$ in this algebra, so that
\begin{equation}\label{deprel}
\Ess f(t) = \sum_{k=1}^N c_k \M_{y_k} \T_{x_k} f(t) = \sum_{k=1}^N
c_k e^{2\pi i y_k \cdot t} f(t-x_k) = 0
\end{equation}
for almost every $t\in \BBR^d$, some constants $0 \neq c_k \in
\BBC$, and some points $(x_k,y_k) \in \Gamma_0$.  Of course, this
means precisely that the Gabor system $\mathcal G \big(f, \{(x_k,
y_k)\}_{1\leq k \leq N}\big)$ is linearly dependent.

In the product setting, we will deduce Theorem
\ref{maintheorem} (and Theorem \ref{lptheorem}) from the following more general result, whose proof
is modeled after a basic argument from the spectral theory of random Schr\"odinger operators
(cf.\ \textit{e.g.\ }Lemma V.2.1 of \cite{carmonalacroix}).

\begin{theorem}\label{generaltheorem}
Suppose $\Gamma$ and $\Lambda$ are arbitrary lattices in
$\BBR^d$.  Let $\gamma_1, \ldots, \gamma_N$ be \emph{distinct}
elements of $\Gamma$, and let $\psi_1, \ldots, \psi_N \in
L^\infty(\BBR^d / \Lambda)$ be nonzero Lebesgue-almost everywhere,
viewed as $\Lambda$-periodic functions on $\BBR^d$.  Then if $0 < p \leq 2$, the
operator $\Ess$ on $L^p(\BBR^d)$ defined by
\[\Ess f(x) = \sum_{k=1}^N \psi_k(x) f(x + \gamma_k)\] has kernel
$\ker(\Ess) = \{0\}$.
\end{theorem}

\begin{remark}\nonumber
Firstly, we note that the relevant \emph{phase space} lattice in the
context of Theorem \ref{maintheorem} is $\Gamma_0 = \Gamma \times
\Lambda^{\perp}$, where $\Lambda^{\perp} \leq \widehat{\BBR^d} \cong
\BBR^d$ is the annihilator or ``dual lattice'' of $\Lambda$.
Secondly, we remark that the generalization from characters to more
general $L^\infty(\BBR^d / \Lambda)$ functions is invited by the
requirement that $\gamma_1, \ldots , \gamma_N$ be distinct, which is
necessary for our proof; notice that if one groups together all
terms associated to a common $x_k$ in (\ref{deprel}), one obtains a
similar expression in which the $\psi_k$ are trigonometric
polynomials (provided $\Gamma_0$ is a product lattice).
\end{remark}

\begin{proof}[Proof of Theorem \ref{generaltheorem}]
For $f \in L^p(\BBR^d)$, we begin by examining $\Ess f$ along
$\Gamma$-orbits; specifically, since $\Gamma \leq \BBR^d$ is discrete, for almost every $x \in \BBR^d$ we
have a sequence $u_x \in \ell^p\Gamma$ defined by $u_x(\gamma) = f(x
+ \gamma)$.  Of course, since $p\leq 2$, we have $\ell^p \Gamma \subseteq \ell^2 \Gamma$; this accounts for the extension to more general $L^p$ spaces in Theorem \ref{lptheorem}.  Then we have $\Ess f(x + \gamma) = \sum_{k=1}^N
\psi_k(x + \gamma) u_x (\gamma + \gamma_k)$ for almost every such
$x$, and accordingly for almost every $x \in \BBR^d$ we can study
the operator $\Ess_x \in \mathcal B (\ell^2\Gamma)$ given by
\[\Ess_x u(\gamma) = \sum_{k=1}^N \psi_k (x+\gamma) u(\gamma +
\gamma_k).\]  If $f \in \ker(\Ess)$, then of course we must have
$u_x \in \ker(\Ess_x)$ for almost every $x$; thus, to prove the
theorem it suffices to show that $\ker(\Ess_x) = \{0\} \subset
\ell^2 \Gamma$ for almost every $x\in \BBR^d$.

Now the key point in requiring $\Lambda$ to be a lattice is that the
family of operators $\Ess_x$ is naturally parametrized by a finite
measure space; namely, it is clear that $\Ess_x$ depends only on the
class of $x$ in $\BBR^d / \Lambda$, so we can consider a Borel
measurable family $x\mapsto \Ess_x$ from $\BBR^d / \Lambda$ to the
space $\mathcal B(\ell^2 \Gamma)$ equipped with the operator norm
topology. Moreover, this family carries a natural action of the
group $\Gamma$. Indeed, for any $\gamma_0 \in \Gamma$, let
$\T_{\gamma_0} \in \mathcal U(\ell^2 \Gamma)$ denote the unitary
translation operator given by $\T_{\gamma_0} u(\gamma) = u(\gamma -
\gamma_0)$. Then one can readily check the commutation relation
\[\T_{\gamma_0} \Ess_x \T_{\gamma_0}^* = \Ess_{x-\gamma_0},\] where
the ``$x-\gamma_0$'' should be interpreted via the obvious
measure-preserving action of $\Gamma \leq \BBR^d$ on $\BBR^d
/ \Lambda$. Since the translations $\T_{\gamma_0}$ are unitary, this
commutation relation descends to kernel projections; that is,
setting $\E(x)$ to be the orthogonal projection onto $\ker(\Ess_x)
\subset \ell^2 \Gamma$, we have
\begin{equation}\label{conjugation}
\T_{\gamma_0} \E(x) \T_{\gamma_0}^* = \E(x-\gamma_0)
\end{equation}
for almost every $x \in \BBR^d/\Lambda$ and all $\gamma_0 \in
\Gamma$.

Again, the goal is to show that $\ker(\Ess_x) = \{0\}$ almost surely
in $x\in \BBR^d / \Lambda$; equivalently, we want to show $\tr \,
\E(x) = \dim \ker(\Ess_x) = 0$ almost everywhere, where ``$\tr$''
denotes the usual trace on $\mathcal B(\ell^2 \Gamma)$. This trace
function is measurable on $\BBR^d / \Lambda$; we defer the proof of
this fact to the next section.  Now for any fixed $\gamma_0 \in
\Gamma$, one has
\[\tr \, \E(x) = \sum_{\gamma \in \Gamma} \left \langle \E(x) \,
\delta_\gamma \, , \, \delta_\gamma \right \rangle = \sum_{\gamma
\in \Gamma} \left \langle \T_{\gamma} \E(x) \T_{\gamma}^* \,
\delta_{\gamma_0} \, , \, \delta_{\gamma_0} \right \rangle =
\sum_{\gamma \in \Gamma} \left \langle \E(x-\gamma) \,
\delta_{\gamma_0} \, , \, \delta_{\gamma_0}\right \rangle,\] where
$\delta_\gamma \in \ell^2 \Gamma$ denotes the characteristic
function of $\{\gamma\}$.  Since the translation action of $\Gamma$
preserves the natural measure on $\BBR^d / \Lambda$, one can
eliminate it from the right-hand side by averaging:
\begin{align*}
\int_{\BBR^d / \Lambda} \tr\, \E(x) \, \mathrm dx &= \sum_{\gamma
\in \Gamma} \int_{\BBR^d / \Lambda} \left \langle \E(x-\gamma) \,
\delta_{\gamma_0} \, , \, \delta_{\gamma_0} \right \rangle \,
\mathrm dx\\ &= \sum_{\gamma \in \Gamma} \int_{\BBR^d / \Lambda}
\left \langle \E(x) \, \delta_{\gamma_0} \, , \,
\delta_{\gamma_0}\right \rangle \, \mathrm dx. \nonumber
\end{align*}
Now of course the summand on the right-hand side is nonnegative and
independent of $\gamma \in \Gamma$, so the sum must be either zero
or infinite.  On the other hand, the integrand on the left-hand side
is nonnegative, so if one could deduce that the sum were actually
zero, one would obtain $\tr \, \E(x) = 0$ for almost every $x \in
\BBR^d / \Lambda$ as desired.  Thus, one would like to deduce \textit{a
priori} that the integral on the left-hand side is \emph{finite},
which could be accomplished by showing $\tr \, \E(x) \lesssim 1$.

For the case $d=1$ this bound is easily achieved as follows.  By
metaplectic reductions we may assume for simplicity that $\Gamma =
\BBZ \leq \BBR$; thus the operators $\Ess_x$ take the form \[\Ess_x
u(n) = \sum_{k=1}^N \psi_k(x + n) u(n + m_k)\] for some $m_1 <
\ldots < m_N \in \BBZ$.  Consider an arbitrary $u \in \ker (\Ess_x)
\subset \ell^2\BBZ$.  Since we assume $\psi_k \neq 0$ almost
everywhere, we may assume that in fact $\psi_k (x+n) \neq 0$ for all
$1 \leq k \leq N$ and all $n \in \BBZ$.  Thus one obtains a
\emph{recurrence relation} for the values of $u$; for any $n \in
\BBZ$, the equation $\Ess_x u(n) = 0$ determines $u(n + m_j)$ in
terms of the values $u(n+m_k)$, $k \neq j$.  In particular, the
function $u \in \ell^2\BBZ$ is completely determined by its values
on the interval $\{n \; | \; m_1 \leq n \leq m_N\} \subset \BBZ$,
and hence \[\tr \, \E(x) = \dim \ker (\Ess_x) \leq m_N - m_1 + 1
\lesssim 1\] for almost every $x$, as desired.  This completes the
proof of Theorem \ref{generaltheorem} for $d=1$.

\smallskip
Unfortunately, however, no such argument yields an \emph{a priori}
estimate $\tr \, \E(x) \lesssim 1$ for $d \geq 2$; this failure
should be apparent from the discussion of the next section.  A
na\"ive remedy for this situation is simply to force uniform
boundedness of the traces by cutting each $\E(x)$ with a common
finite-rank projection.  To this end, for $A \subset \Gamma$ let
$\chi_A$ denote the characteristic function of $A$, and define the
operator $\E_A(x) := \chi_A \E(x)$; here $\chi_A \in \ell^\infty
\Gamma$ is viewed as a projection operator on $\ell^2 \Gamma$. Again
using the commutation relation (\ref{conjugation}), for any
$\gamma_0 \in \Gamma$ we have
\[\tr \, \E_A(x) = \sum_{\gamma \in A} \left \langle \E(x) \, \delta_\gamma \, , \, \delta_\gamma \right \rangle = \sum_{\gamma \in A - \gamma_0} \left \langle \E(x + \gamma) \, \delta_{\gamma_0} \, , \, \delta_{\gamma_0} \right \rangle.\]

Just as above, we integrate over $\BBR^d / \Lambda$ to obtain
\[|A-\gamma_0| \int_{\BBR^d / \Lambda} \left \langle \E(x) \,
\delta_{\gamma_0} \, , \, \delta_{\gamma_0}\right \rangle \, \mathrm
dx = \int_{\BBR^d / \Lambda} \tr \, \E_A(x)\, \mathrm dx \lesssim
\|\tr \, \E_A\|_{L^\infty (\BBR^d / \Lambda)},\] whence

\begin{equation}\label{tracebound}
\int_{\BBR^d / \Lambda} \left \langle \E(x) \,
\delta_{\gamma_0} \, , \, \delta_{\gamma_0} \right \rangle \,
\mathrm dx \lesssim \frac{\|\tr \, \E_A\|_\infty}{|A|}
\end{equation}
 with the
implied constant independent of $A$.  (Here ``$|A|$'' denotes the
cardinality of $A \subset \Gamma$.)  By the crucial ``growth''
property of $\Gamma$ given by Lemma \ref{propagationlemma} of the
next section, we can find a sequence of sets $A_n \subset \Gamma$
such that
\begin{equation}\label{dimensiondecay}
\frac{|\tr \, \E_{A_n}(x)|}{|A_n|} \longrightarrow 0
\end{equation}
uniformly in almost every $x$ as $n \rightarrow \infty$.  The basic idea is that
while it is not immediately obvious that $\dim \ker(\Ess_x) <
\infty$, we can still determine a \emph{large} number of values of
$u \in \ker(\Ess_x)$ from a relatively \emph{small} number of its
values.  Exploiting these sets $A_n$ in the estimate
(\ref{tracebound}), we obtain $\langle \E(x) \, \delta_{\gamma_0} \,
, \, \delta_{\gamma_0}\rangle = 0$ for almost every $x$.  Since
$\gamma_0 \in \Gamma$ was arbitrary, we have
\[\dim \ker(\Ess_x) = \tr \, \E(x) = \sum_{\gamma \in \Gamma} \left
\langle \E(x) \, \delta_\gamma \, , \, \delta_\gamma \right \rangle
= 0\] for almost every $x \in \BBR^d / \Lambda$; modulo the two
deferred claims, this completes the proof of the theorem.
\end{proof}

\section{Proofs of auxiliary results}
It remains to verify the measurability of the function $\tr \, \E:
\BBR^d /\Lambda \rightarrow \BBC$ and to produce a sequence of sets
$A_n \subset \Gamma$ satisfying the property (\ref{dimensiondecay}).
The first task should be viewed as a technicality and may well be
standard fare; we were unable to find the specific result we require in the literature, but see for example Chapter 5 of \cite{carmonalacroix} for a discussion of related measurability issues.
On the other hand, the second task seems more
fundamental to Theorem \ref{generaltheorem}.  Throughout this
section, we employ the same notations as those used in the proof
above.

\subsection{Measurability of the kernel projection trace}
The idea behind proving measurability of $\tr \, \E$ is simply to
use functional calculus for the self-adjoint operators $\Ess_x^*
\Ess_x$, since $\ker(\Ess_x) = \ker(\Ess_x^*\Ess_x)$, so that $\E(x)
= \operatorname{proj}_{\ker(\Ess_x)}$ is a spectral projection of
$\Ess_x^* \Ess_x$.  We begin with the following general lemma:

\begin{lemma}\label{spectralmeasurelemma}
Let $\mathcal H$ be a separable Hilbert space, and fix $u,v \in
\mathcal H$.  For an operator $\Ess \in \mathcal B(\mathcal H)$, let
$\mu_\Ess$ denote the spectral measure on $\BBR$ associated to the
self-adjoint operator $\Ess^* \Ess$, $u$, and $v$, so that
\[\left \langle (\Ess^* \Ess)^k u, v \right \rangle = \int_\BBR x^k \, \mathrm d
\mu_\Ess (x).\] Then the map $\Ess \mapsto \mu_\Ess$ is continuous
from $\mathcal B(\mathcal H)$ equipped with the strong-$*$ operator
topology to the space of measures on $\BBR$ equipped with the
weak-$*$ topology.
\end{lemma}

\begin{proof}
Suppose we have a sequence $\Ess_n \in \mathcal B(\mathcal H)$
converging to $\Ess \in \mathcal B(\mathcal H)$ in the strong-$*$
operator topology, so that $\Ess_n \xi \rightarrow \Ess \xi$ and
$\Ess_n^* \xi \rightarrow \Ess^* \xi$ for all $\xi \in \mathcal H$.
By an application of the uniform boundedness principle, note that
all the spectral measures $\mu_{\Ess_n}$ and $\mu_\Ess$ are
supported on a common compact interval $I \subset \BBR$.  Now for
any polynomial $p\in \BBC[x]$, the strong-$*$ convergence implies
\textit{a fortiori} that $p(\Ess_n^* \Ess_n) \rightarrow
p(\Ess^*\Ess)$ in the weak operator topology on $\mathcal B(\mathcal
H)$; in particular, we have \[\int_I p(x) \, \mathrm d
\mu_{\Ess_n}(x) \longrightarrow \int_I p(x) \, \mathrm d \mu_\Ess
(x).\]  An application of the Weierstrass polynomial approximation
theorem now shows that $\mu_{\Ess_n} \rightarrow \mu_\Ess$ in the
weak-$*$ topology.
\end{proof}

With this lemma in hand, we can now check the measurability of $\tr
\, \E$.  By the definition of the operators $\Ess_x$ and the
measurability of the functions $\psi_k$ appearing therein, it is
clear that for each $u \in \ell^2 \Gamma$ the quantities $\|\Ess_x
u\|_2$ and $\|\Ess_x^*u\|_2$ are measurable functions of $x$; thus,
the map $x \mapsto \Ess_x$ is Borel measurable from $\BBR^d/\Lambda$
to $\mathcal B(\ell^2\Gamma)$ equipped with the strong-$*$ operator
topology.  (It is easy to check that this implies the map is still
Borel when $\mathcal B(\ell^2\Gamma)$ is equipped with the operator
norm topology, as claimed in the previous section.)  Now fixing $u,v
\in \ell^2\Gamma$ and invoking Lemma \ref{spectralmeasurelemma}, the
map $x \mapsto \int \varphi \, \mathrm d\mu_{\Ess_x}$ is measurable
for all $\varphi \in C_c(\BBR)$; thus, approximating the
characteristic function $\chi_{\{0\}}$ as a pointwise limit of
continuous functions, so is the map \[x \longmapsto
\mu_{\Ess_x}(\{0\}) = \langle \proj_{\ker(\Ess_x)} u,v \rangle =
\langle \E(x) u,v\rangle.\]  In other words, the projection-valued
map $\E$ is weakly measurable, whence it trivially follows that $\tr
\, \E: \BBR^d / \Lambda \rightarrow \BBC$ is a Borel measurable
function.

\subsection{Sets $A_n \subset \Gamma$ with property
(\ref{dimensiondecay})}  To complete the proof of Theorem \ref{generaltheorem}, we need to find
a sequence of subsets $A_n \subset \Gamma$ satisfying
(\ref{dimensiondecay}), namely \[\frac{|\tr \, \E_{A_n}(x)|}{|A_n|}
\longrightarrow 0\] as $n\rightarrow \infty$, uniformly in almost
every $x$.

For any finite $A \subset \Gamma$, let $\mathcal K_A \subset \ell^2 \Gamma$ denote the subspace of functions supported on $A$; thus $\mathcal K_A \cong \BBC^{|A|}$.  Since $\E_A(x) = \chi_A \E(x)$, the range of $\E_A(x)$ is a subspace of $\mathcal K_A$, and the restriction $\E_A (x) \chi_A = \chi_A \E(x) \chi_A$ of $\E_A (x)$ to $\mathcal K_A$ is a self-adjoint operator of norm at most $1$; hence \[\tr \, \E_A(x) = \tr \, \big(\E_A(x) \chi_A\big) \leq \dim \, \operatorname{ran} \, \big( \E_A(x) \chi_A\big) \leq \dim \big(\chi_A \ker (\Ess_x)\big).\]  Thus it suffices to find $A_n \subset \Gamma$ with \[\frac{\dim \big(\chi_{A_n} \ker (\Ess_x)\big)}{|A_n|} \longrightarrow 0.\]

To produce the sets $A_n$, we take the same basic perspective
used to show $\tr \, \E(x) \lesssim 1$ in the $d=1$ setting:  For any
$u \in \ker(\Ess_x)$ and $\gamma \in \Gamma$, the definition of
$\Ess_x$ allows one automatically to determine $u(\gamma +
\gamma_k)$, provided one knows $u(\gamma + \gamma_j)$ for all $j
\neq k$, $1 \leq j \leq N$.\footnote{The uniformity in $x$ of the
decay follows trivially from the proof below and will not be
mentioned; just as in the $d=1$ case treated above, the only
implicit mention of $x$ in the argument is the requirement that
$\psi_k(x+n) \neq 0$ for all $k$ and $n$, which as before can be
guaranteed by eliminating a measure-zero set of $x$.}

To rephrase this perspective slightly, fix any subset $C_0 \subset
\Gamma$ and any element $\gamma_0 \in C_0$.  Then for any subset $C
\subset \Gamma$, we define $\mathcal P_{C_0, \gamma_0} (C)$ to be
the minimal subset of $\Gamma$ satisfying $C \subset \mathcal
P_{C_0, \gamma_0}(C)$ and the implication \[\gamma + \left(C_0 \setminus
\{\gamma_0\}\right) \subset \mathcal P_{C_0,\gamma_0}(C) \,
\Longrightarrow \gamma + \gamma_0 \in \mathcal P_{C_0,
\gamma_0}(C).\]  In other words, the larger set $\mathcal
P_{C_0,\gamma_0}(C)$ is iteratively ``grown from $C$'' by the
following rule:  ``At any given stage, if the set contains a
translate of $C_0 \setminus \{\gamma_0\}$, put that same translate
of $\gamma_0$ in the set and continue.''  Now by the above
reasoning, if we take $C_0 = \{\gamma_1, \ldots, \gamma_N\}$ and set
$\gamma_0 = \gamma_k$ for some $1\leq k \leq N$, we see that
\[\dim \big(\chi_{\mathcal P_{C_0,
\gamma_k}(C)} \ker(\Ess_x) \big) \leq \dim \big(\chi_C \ker(\Ess_x) \big) \leq |C|\] for any finite $C \subset \Gamma$.  The
following property of $\Gamma$ will thus give the desired decay
(\ref{dimensiondecay}).

\begin{lemma}\label{propagationlemma}
Let $\Gamma \leq \BBR^d$ be a lattice, and fix an arbitrary finite
subset $C_0 \subset \Gamma$.  Then there exists a $\gamma_0 \in C_0$
and a sequence of subsets $C_n \subset \Gamma$, $n \in \BBN$, such
that $|C_n| \lesssim n^{d-1}$ and $|\mathcal P_{C_0, \gamma_0}(C_n)|
\gtrsim n^d$, with $\mathcal P_{C_0, \gamma_0}$ defined as above and
the implied constants depending only on $d$, $\Gamma$, and $C_0
\subset \Gamma$.
\end{lemma}

The following proof is essentially due to Thiele (\cite{thielecomm})
and is best understood geometrically; the reader may find it rather
helpful to draw some pictures for the case $\Gamma = \BBZ \times
\BBZ \leq \BBR^2$, as doing so will highlight the simplicity of the
argument.

\begin{proof}
First we note that the conclusion of the lemma is clearly invariant
under translations of the set $C_0$, possibly with the exception of
the implied constants in the cardinality estimates.  It will be
apparent from the proof below that these constants are indeed
unaffected by translations of $C_0$, so for simplicity we may assume
$0 \in C_0$.

We will choose $\gamma_0 \in C_0$ to be an extreme point of the
convex hull $\ch (C_0)$ of $C_0 \subset \BBR^d$; again by
translation-invariance of the claim, we may assume $\gamma_0 = 0 \in
\BBR^d$.  Then it is a matter of routine to check that there is a
rank-$(d-1)$ subgroup $K \leq \Gamma$ with $K \cap C_0 = \{0\}$;
geometrically, one should view $K$ as the intersection of a
hyperplane in $\BBR^d$ with the lattice $\Gamma$.  The quotient
$\Gamma / K$ is thus cyclic, and we can choose $x \in \Gamma$ such
that $x + K$ generates $\Gamma / K$ and such that \[C_0 \setminus
\{0\} \subset (x + K) \cup (2x + K) \cup \ldots \cup (mx + K) =:
\bar C\] for some $m \in \BBN$.  This latter set $\bar C$ should be
viewed as a stack of hyperplanes in $\Gamma$ adjacent to $K$ that
foliate a certain ``strip'' in the lattice $\Gamma$.

The essential observation is simply that the definition of $\mathcal
P_{C_0, \gamma_0}$ immediately yields $\mathcal P_{C_0,0}(\bar C)
\supset \bar C \cup K$; in particular, $\mathcal P_{C_0,0}(\bar C)$
contains the shifted hyperplane stack $-x + \bar C$.  By induction,
we obtain \[\mathcal P_{C_0,0} (\bar C) \supset
\bigcup_{j=-\infty}^0 (jx + K) =: H,\] which is the intersection of
$\Gamma$ with a half-space in $\BBR^d$.  In short, an entire
half-space in $\Gamma$ can be grown from the codimension-$1$
``strip'' $\bar C$ by the procedure used to define $\mathcal
P_{C_0,0}(\bar C)$. Passing from this observation to the
quantitative statement of the lemma is a simple matter, along the following lines.

Let $B_r(0)$ denote the usual ball of radius $r>0$ centered at $0$
in $\BBR^d$, set $C_n := B_n(0) \cap \bar C$, and let $\delta$ denote the (Euclidean) diameter of $C_0$.  Using the iterative definition of $\mathcal P_{C_0, \gamma_0}$, one checks that \[\mathcal P_{C_0,0}(C_n) \supset jx + \big(B_{n + (j-1)\cdot 10\delta}(0) \cap K\big)\] for each $j\in \BBZ$ satisfying $-\frac{n}{100\delta} < j \leq 0$, say.  In particular, we have \[\mathcal P_{C_0,0}(C_n) \supset \bigcup_{-\frac{n}{200\delta} + 1 \leq j \leq -1} \Big(jx + \big(B_{\frac{n}{2}}(0)\cap K\big)\Big)\] for all sufficiently large $n$.  Thus, since $K$ is a rank-$(d-1)$ subgroup of the lattice $\Gamma \leq \BBR^d$, we have $|C_n|
\lesssim_{d,\Gamma,C_0} n^{d-1}$ and $|\mathcal P_{C_0,0} (C_n)|
\gtrsim_{d,\Gamma,C_0} n^d$, as desired.
\end{proof}

Finally, of course, we can set $C_0 = \{\gamma_1 , \ldots ,
\gamma_N\}$, apply Lemma \ref{propagationlemma}, and take $A_n$ to
be the resulting sets $A_n = \mathcal P_{C_0,\gamma_0}(C_n)$ to
obtain property (\ref{dimensiondecay}); this completes the proof of
Theorem \ref{generaltheorem} and hence that of Theorem
\ref{maintheorem}.


\section{Additional remarks} \label{remarks}
\begin{remark}\label{covolumeoneremark}
It is somewhat interesting to note how the proof of Theorem \ref{generaltheorem} degenerates for the special case of $\Lambda = \Gamma$, which in the $d=1$ setting corresponds to that of covolume-$1$ phase space lattices in $\BBR \times \widehat \BBR$ (after exploiting metaplectic symmetries).  In this case, the action of $\Gamma$ on $\BBR^d / \Gamma$ is of course trivial, so each $\Ess_x$ is simply a linear combination of translation operators on $\ell^2 \BBZ$.  One need not consider the kernel projections $\E (x)$, and modulo an application of the Fourier transform our proof becomes essentially that of Proposition 2 in \cite{hrtpaper}.
\end{remark}

\begin{remark}\label{productremark}
While the particular structure of product lattices was clearly exploited in the proof of Theorem \ref{maintheorem}, the difficulties posed by general lattices for our method of proof may still not be apparent.  Indeed, by a rather clever combination of metaplectic transformations and passage to higher dimensions, in \cite{linnell99} Linnell shows that one can reduce considerations to lattices $\Gamma_0 \leq \BBR^d \times \BBR^d$ such that $\Gamma_0 \cap \{0\} \times \BBR^d = \{0\} \times \BBZ^d$; such lattices are
temptingly close to product lattices, as their ``translation
components'' (\textit{i.e.\ }their projections to $\BBR^d \times
\{0\}$) are lattices in $\BBR^d$. However, even these lattices seem
in general to be out of reach for our techniques. As a particular
example, in the case $d=2$, consider a lattice $\Gamma_0 \leq
\BBR^4$ with $\BBZ$-basis
\[\{(0,0,1,0)\, , \, (0,0,0,1) \, , \, (\alpha_1,\alpha_2,\alpha_3,\alpha_4) \, , \,
(\beta_1,\beta_2,\beta_3,\beta_4)\},\] where $\alpha_1$, $\alpha_2$, $\beta_1$, $\beta_2$, and $\alpha_1 \beta_3 + \alpha_2 \beta_4 - \alpha_3 \beta_1 - \alpha_4 \beta_2 \in \BBR$ are linearly independent over $\BBZ$.   Such a lattice is of the form to which Linnell
reduces; on the other hand, by checking the values of the symplectic
form on pairs of basis elements, one sees that such a $\Gamma_0$ cannot be
symplectically equivalent to a product lattice.  (A product
lattice in $\BBR^4$ must have a $\BBZ$-basis $\{v_1,v_2,v_3,v_4\}$ with
$[v_1,v_2] = [v_3,v_4] = 0$, where ``$[\cdot,\cdot]$'' denotes the
symplectic form on $\BBR^4$.  Routine algebra shows that our choice of $\Gamma_0$ has no such basis; we thank Nets Katz for pointing out this
efficient method of generating examples of lattices that are not
symplectically equivalent to product lattices.)  Moreover, for such a $\Gamma_0$ with $(\alpha_3 , \alpha_4)$ and $(\beta_3, \beta_4)$ generating an infinite-covolume subgroup of $\BBR^2$,
the reader is invited to carry out a similar analysis to that in the
arguments above and see what goes wrong. In short, it is not clear
to us how to parametrize the operators $\Ess_x$ arising in the
proof by a \emph{finite} measure space that admits a
measure-preserving action of the relevant translation lattice.
\end{remark}

\begin{remark}\label{groupremark}
As a further comment on the nature of the product lattice setting, we note that the ``post-metaplectic reduction'' Theorem
\ref{generaltheorem} can easily be reformulated for pairs of
lattices in more general locally compact (not necessarily abelian)
groups. The only part of the proof that causes any difficulty in
extension is the analogue of Lemma \ref{propagationlemma} for the
\emph{translation} lattice; we hope to address generalizations of
these results in future work. In this light, modulo metaplectic
symmetries, Theorem \ref{maintheorem} could be viewed as mainly
being a theorem ``about the translation lattice''; thus, its role as
evidence for the full HRT Conjecture might be somewhat dubious.
\end{remark}

\bibliography{latticebib}{}

\begin{thebibliography}{HRT96}

\bibitem[BS09]{bownikspeegle}
Marcin Bownik and Darrin Speegle.
\newblock Linear independence of {P}arseval wavelets.
\newblock Preprint, 2009.

\bibitem[CL90]{carmonalacroix}
Ren{\'e} Carmona and Jean Lacroix.
\newblock {\em Spectral theory of random {S}chr\"odinger operators}.
\newblock Probability and its Applications. Birkh\"auser Boston Inc., Boston,
  MA, 1990.

\bibitem[Dem10]{demetermn}
Ciprian Demeter.
\newblock Linear independence of time frequency translates for special
  configurations.
\newblock {\em Math. Res. Lett.}, 17(4):761--779, 2010.

\bibitem[DZ10]{demeterzaharescu}
Ciprian Demeter and Alexandru Zaharescu.
\newblock Proof of the {H}{R}{T} conjecture for (2,2) configurations.
\newblock Preprint, arxiv:1006.0735 [math.CA], 2010.

\bibitem[Hei06]{heil}
Christopher Heil.
\newblock Linear independence of finite {G}abor systems.
\newblock In {\em Harmonic analysis and applications}, Appl. Numer. Harmon.
  Anal., pages 171--206. Birkh\"auser Boston, Boston, MA, 2006.

\bibitem[HRT96]{hrtpaper}
Christopher Heil, Jayakumar Ramanathan, and Pankaj Topiwala.
\newblock Linear independence of time-frequency translates.
\newblock {\em Proc. Amer. Math. Soc.}, 124(9):2787--2795, 1996.

\bibitem[Jit99]{jitomirskaya}
Svetlana~Ya. Jitomirskaya.
\newblock Metal-insulator transition for the almost {M}athieu operator.
\newblock {\em Ann. of Math. (2)}, 150(3):1159--1175, 1999.

\bibitem[Kut02]{kutyniok}
Gitta Kutyniok.
\newblock Linear independence of time-frequency shifts under a generalized
  {S}chr\"odinger representation.
\newblock {\em Arch. Math. (Basel)}, 78(2):135--144, 2002.

\bibitem[Lin91]{linnell91}
Peter~A. Linnell.
\newblock Zero divisors and group von {N}eumann algebras.
\newblock {\em Pacific J. Math.}, 149(2):349--363, 1991.

\bibitem[Lin99]{linnell99}
Peter~A. Linnell.
\newblock von {N}eumann algebras and linear independence of translates.
\newblock {\em Proc. Amer. Math. Soc.}, 127(11):3269--3277, 1999.

\bibitem[Ste93]{bigstein}
Elias~M. Stein.
\newblock {\em Harmonic analysis: real-variable methods, orthogonality, and
  oscillatory integrals}, volume~43 of {\em Princeton Mathematical Series}.
\newblock Princeton University Press, Princeton, NJ, 1993.
\newblock With the assistance of Timothy S. Murphy, Monographs in Harmonic
  Analysis, III.

\bibitem[Thi07]{thielecomm}
Christoph Thiele, 2007.
\newblock Personal communication.

\end{thebibliography}
\bibliographystyle{alpha}
\end{document}